\documentclass[10pt,final]{amsart}
\usepackage{amsmath,amssymb,amsthm,showkeys}

\usepackage{amsmath,amssymb,amsthm,showkeys}

\usepackage{amsmath}
\usepackage{graphicx}
\usepackage{color}

\begin{document}

\newtheorem{lem}{Lemma}[section]
\newtheorem{prop}{Proposition}[section]
\newtheorem{cor}{Corollary}[section]
\numberwithin{equation}{section}
\newtheorem{thm}{Theorem}[section]

\newtheorem{ass}{Assumption}[section]

\theoremstyle{remark}
\newtheorem{example}{Example}[section]
\newtheorem*{ack}{Acknowledgments}

\theoremstyle{definition}
\newtheorem{definition}{Definition}[section]

\theoremstyle{remark}
\newtheorem*{notation}{Notation}
\theoremstyle{remark}
\newtheorem{remark}{Remark}[section]

\newenvironment{Abstract}
{\begin{center}\textbf{\footnotesize{Abstract}}%
\end{center} \begin{quote}\begin{footnotesize}}
{\end{footnotesize}\end{quote}\bigskip}

\newcommand{\triple}[1]{{|\!|\!|#1|\!|\!|}}

\newcommand{\xx}{\langle x\rangle}
\newcommand{\ep}{\varepsilon}
\newcommand{\al}{\alpha}
\newcommand{\be}{\beta}
\newcommand{\de}{\partial}
\newcommand{\la}{\lambda}
\newcommand{\La}{\Lambda}
\newcommand{\ga}{\gamma}
\newcommand{\del}{\delta}
\newcommand{\Del}{\Delta}
\newcommand{\sig}{\sigma}
\newcommand{\ome}{\omega}
\newcommand{\Ome}{\Omega}
\newcommand{\C}{{\mathbb C}}
\newcommand{\N}{{\mathbb N}}
\newcommand{\Z}{{\mathbb Z}}
\newcommand{\R}{{\mathbb R}}
\newcommand{\Rn}{{\mathbb R}^{n}}
\newcommand{\Rnu}{{\mathbb R}^{n+1}_{+}}
\newcommand{\Cn}{{\mathbb C}^{n}}
\newcommand{\spt}{\,\mathrm{supp}\,}
\newcommand{\Lin}{\mathcal{L}}
\newcommand{\SSS}{\mathcal{S}}
\newcommand{\F}{\mathcal{F}}
\newcommand{\xxi}{\langle\xi\rangle}
\newcommand{\eei}{\langle\eta\rangle}
\newcommand{\xei}{\langle\xi-\eta\rangle}
\newcommand{\yy}{\langle y\rangle}
\newcommand{\dint}{\int\!\!\int}
\newcommand{\hatp}{\widehat\psi}

\renewcommand{\Re}{\;\mathrm{Re}\;}
\renewcommand{\Im}{\;\mathrm{Im}\;}

\title[Scale invariant  energy smoothing  estimates]{Scale  invariant  energy smoothing estimates
 for the  Schr\"odinger Equation with small Magnetic
Potential}

\author{Vladimir Georgiev }
\address{ Vladimir Georgiev,  Dipartimento di Matematica, Universit\`a degli Studi di Pisa,
Largo Bruno Pontecorvo 5, 56127 Pisa, Italy.}


\author{Mirko Tarulli }
 \address{ Mirko Tarulli,  Dipartimento di Matematica, Universit\`a degli Studi di Pisa,
Largo Bruno Pontecorvo 5, 56127 Pisa, Italy.}

\thanks{The  authors are
  partially
supported by Research Training Network (RTN) HYKE,  financed by
the European Union, contract number : HPRN-CT-2002-00282.}

\keywords{AMS Subject Classification: %
 Key words: Schr\"odnger equation, smoothing estimates, magnetic
potential}

\begin{abstract}
We consider  some scale invariant generalizations of the smoothing
estimates for the free Schr\"odnger equation obtained by Kenig,
Ponce and Vega in \cite{KPV93}, \cite{KPV98}. Applying these
estimates and using appropriate commutator estimates, we obtain
similar scale invariant smoothing estimates for perturbed
Schr\"odnger equation with small magnetic potential.
\end{abstract}

\maketitle

\date{}



\section{Introduction}

In this work we study smoothing properties of the Schr\"odinger
equation with magnetic potential
 $$A = (A_1(t,x),\cdots,A_n(t,x)), \ \ x\in \mathbb{R}^n.$$
 Here $A_j(t,x), j=1,\cdots,n,$ are real
valued functions, $n \geq 3$ and the corresponding Cauchy problem
for  Schr\"odinger equation has the form
\begin{equation}\label{eq.lwa}
\begin{cases}
\partial_t u - i\Delta_A u = F,\ \ \  t \in \mathbb{R} , \ \ \ x \in \mathbb{R}^n\\
u(0,x)=f(x),
\end{cases}
\end{equation}
where
\begin{equation}\label{eq.covl}
    \Delta_A
    = \sum_{j=1}^n (\partial_{x_j} - i A_j) (\partial_{x_j} - i
    A_j).
\end{equation}

The energy type estimates and well -- posedness of the Cauchy
problem \eqref{eq.lwa} in energy space are studied in the works
\cite{Doi94} and \cite{Doi96} of Doi.

Since the smoothing properties of this evolution problem are
closely connected with suitable resolvent estimates for the
solution $U=U(x)$ of the elliptic problem
\begin{equation}\label{eq.res1}
\begin{cases}
\varepsilon U - i\Delta_A U - i\tau U = H,\ \ \  \varepsilon >0,
\tau >0, \ \ \ x \in \mathbb{R}^n,  H=H(x),
\end{cases}
\end{equation}
we can use as a starting point the scale invariant smoothing
estimate obtained in the works of Kenig, Ponce, Vega \cite{KPV93}
and Pertham, Vega \cite{KPV98}. This estimate extends earlier
works of Agmon, H\"ormander \cite{AH} and P. Constantin and J.-C.
Saut \cite{CS}.

The scale invariant estimate for \eqref{eq.res1} with $A=0$ has
the form
\begin{equation}\label{eq.KPVsm}
    \|| \nabla_x U \|| \leq C N(H),
\end{equation}
where $C>0$ is independent of $\varepsilon >0, \tau >0,$
$$ \|| G \||^2 = \sup_{R>0} \frac{1}{R} \int_{|x| \leq R} |G(y)|^2
dy$$ is the Morrey - Campanato norm, while
$$
N(H) = \sum_{k \in \mathbb{Z}}2^{k/2} \| H \|_{L^2(2^{k-1} \leq
|x| \leq 2^{k+1})}.$$

From this estimate one can use the simple estimate
$$ \sup_{k \in \mathbb{Z}}2^{-k/2} \| G \|_{L^2(2^{k-1} \leq
|x| \leq 2^{k+1})} \leq C \|| G \||$$ and derive the following
smoothing scale invariant estimate for the solution $u(t,x) $ to
\eqref{eq.lwa} with $A=0$ and $f=0$
\begin{eqnarray}\label{eq.sm}
\int_{\mathbb{R}}\left(  \sup_{k \in \mathbb{Z}} 2^{-k/2}
\|\nabla_x u(t, \cdot) \|_{L^2(2^{k-1} \leq |x| \leq 2^{k+1})}
\right)^2dt \leq \\
\leq C \int_{\mathbb{R}}\left( \sum_{k \in \mathbb{Z}} 2^{k/2} \|
F(t, \cdot) \|_{L^2(2^{k-1} \leq |x| \leq 2^{k+1})} \right)^2 dt.
\nonumber
\end{eqnarray}
Our purpose in this work is to derive similar scale invariant
smoothing estimates for the case of magnetic potential imposing
scale invariant smallness assumptions on the magnetic potential
$A(x).$ The reason why we treat only small magnetic potential is
connected with the necessity to avoid resonances phenomena ( see
\cite{CuSCH}, \cite{DaFa}, \cite{Ku}, \cite{YN}, \cite{Ra},
\cite{Sm}, \cite{Ste04}, \cite{Ta}). The absence of eigenvalues of
$\Delta_A$ with magnetic potential decaying as
$(1+|x|)^{-1-\delta}$ is discussed in  \cite{BMP}. However, even
the remarkable result in \cite{BMP} can not guarantee that $0$ is
not an eigenvalue of the Hamiltonian $\Delta_A.$ The result in
\cite{Hi02} shows that even nontrivial smooth compactly supported
magnetic field can create resonances.

To avoid possible eigenvalues or resonances of $\Delta_A$ we
impose the following assumption on the potential $A.$
\begin{ass} \label{ass2}
There exists $\varepsilon >0 ,$ such that  we have
\begin{equation}\label{eq.ass2}
\max_{1 \leq j \leq n}  \ \sum_{k \in \mathbb{Z}} \sum_{|\beta|
\leq 1} \ 2^{k(1+|\beta|)} \ \|D_x^\beta A_j(t,x)\|_{L^\infty_t
L^\infty_{\{|x| \sim 2^k\}}}\leq \varepsilon.
\end{equation}
\end{ass}

Our main smoothing estimate is the following one.

\begin{thm}
\label{mainsm} There exists $\varepsilon >0$ so that for any
potential $A(x)$ satisfying \eqref{eq.ass2} there exists $C>0,$ so
that  for any $f \in S( \mathbb{R}^n)$ and any $F(t,x) \in
C_0^\infty( \mathbb{R} \times (\mathbb{R}^n  \setminus 0))$ the
solution $u(t,x)$ to \eqref{eq.lwa} satisfies the estimate
\begin{eqnarray}\label{eq.strsmoo}
\int_{\mathbb{R}}\left(  \sup_{k \in \mathbb{Z}} \||x|_k^{-1/2}
u(t,  \cdot) \|_{\dot{H}^{1/2}_x} \right)^2dt &\leq& C \| f \|^2_{
L^{2}_x} +
\nonumber\\
&+& C \int_{\mathbb{R}}\left( \sum_{k \in \mathbb{Z}}
\||x|_k^{1/2}  F(t,  \cdot) \|_{\dot{H}^{-1/2}_x} \right)^2 dt,
\end{eqnarray}
where  $ \dot{H}^{s}_x = \dot{H}^{s}(\mathbb{R}^n )$ is the
classical homogeneous Sobolev space and $|x|^{\pm 1/2}_k =
|x|^{\pm 1/2}Q_k(x)$ and the Paley - Littlewood partition of unity
\begin{equation}\label{eq.PLspace}
    1 = \sum_{k \in \mathbb{Z}} Q_k(x) ,
\end{equation}
is defined as follows
$$
  Q_k(x) =
\varphi \left( \frac{|x|}{2^k}\right),$$ where $\varphi(s) \in
C_0^\infty( (1/2, 2))$ is a non - negative function.
\end{thm}

The key point to derive this estimate is a suitable scale and time
invariant smoothing estimate for the free Schr\"odinger equation

\begin{equation}\label{eq.lwafree}
\begin{cases}
\partial_t u - i\Delta u = F,\ \ \  t >0 , \ \ \ x \in \mathbb{R}^n,\\
u(0,x)=f(x).
\end{cases}
\end{equation}
To be more precise, we introduce the following norms motivated by
the statement of the main result in Theorem \ref{mainsm}. Take
$$ Y= L^2_t( \ell_x^{1,1/2} \dot{H}^{-1/2} ),
\ \ Y^\prime = L^2_t( \ell_x^{\infty,-1/2} \dot{H}^{1/2} ),$$
where the spaces $\ell_x^{q,\alpha}B$ for any Banach space $B$ is
introduced in Section \ref{locdef}. Note that $Y$ is not reflexive
( $(\ell_x^{1,1/2})^\prime = \ell_x^{\infty,-1/2},$ but
$(\ell_x^{\infty,-1/2})^\prime \neq \ell_x^{1,1/2}$).

Then the estimate of the previous theorem can be rewritten in the
form
\begin{equation}\label{eq.lwaeq}
\| u \|^2_{Y^\prime} \leq C \|f\|^2_{L^2} +\| F \|^2_{Y},
\end{equation}
where here and below  \begin{equation}\label{eq.defYYp} Y= L^2_t(
\ell_x^{1,1/2} \dot{H}^{-1/2} ), \ \ Y^\prime = L^2_t(
\ell_x^{\infty,-1/2} \dot{H}^{1/2} ).
\end{equation} We shall call $Y^\prime$
smoothing space.

 Then the main point in the proof of
Theorem \ref{mainsm} is to establish first the following energy
smoothing estimate  for the case $A=0.$
\begin{thm}\label{mainsmfree}
There exists $C>0,$ such that for any $f \in S( \mathbb{R}^n)$ and
any $F(t,x) \in C_0^\infty( \mathbb{R} \times (\R^n \setminus0))$
the solution $u(t,x)$ to \eqref{eq.lwafree} satisfies the estimate
\begin{eqnarray}\label{eq.strsmoofree}
& &\|u\|_{L^\infty_t L^2_x} +\|u\|_{Y^\prime} \leq C \| f \|_{
L^{2}_x} + C \left(\min_{F=F_1+F_2}\|F_1\|_{Y} + \|F_2\|_{L^1_t
L^2_x} \right).
\end{eqnarray}
\end{thm}

On the basis of the estimate in Theorem \ref{mainsmfree} we shall
derive a slightly stronger estimate for the perturbed
Schr\"odinger equation.
\begin{cor}
\label{mainsmmag} There exists $\varepsilon >0$ so that for any
potential $A(x)$ satisfying \eqref{eq.ass2} there exists $C>0,$ so
that  for any $f \in S( \mathbb{R}^n)$ and any $F(t,x) \in
C_0^\infty( \mathbb{R} \times (\mathbb{R}^n  \setminus 0))$ the
solution $u(t,x)$ to \eqref{eq.lwa} satisfies the estimate
\begin{eqnarray}\label{eq.strsmoomag}
& &\|u\|_{L^\infty_t L^2_x} +\|u\|_{Y^\prime} \leq C \| f \|_{
L^{2}_x} + C \left(\min_{F=F_1+F_2}\|F_1\|_{Y} + \|F_2\|_{L^1_t
L^2_x} \right).
\end{eqnarray}
\end{cor}

As an application we consider the following semilinear
Schr\"odinger equation
\begin{equation}\label{eq.lwanon}
\begin{cases}
\partial_t u - i\Delta_A u = |V(t,x) u|^{p},\ \ \  t \in \mathbb{R} , \ \ \ x \in \mathbb{R}^n\\
u(0,x)=f(x),
\end{cases}
\end{equation}
where $p >1$ and $V(t,x)$ is a measurable function satisfying the
inequality
\begin{equation}\label{eq.assV}
      \ \sum_{k \in \mathbb{Z}}  \ 2^{ka} \ \|V(t,x)\|_{L^\infty_t
L^\infty_{\{|x| \sim 2^k\}}}\leq C < \infty.
\end{equation}

Then we have the following global existence result with initial
data having small $L^2-$ norm only.

\begin{thm}
\label{mainsmnonl} Suppose the potential $A(x)$ satisfies
\eqref{eq.ass2},  $V$ obeys \eqref{eq.assV} with $a \in [1,2)$ and
\begin{equation}\label{eq.Aass1}
p = \frac{n+4}{n+2a}.
\end{equation}
 Then there exists $\delta>0,$ so that for any $f \in L^2(
\mathbb{R}^n)$ with
$$ \|f\|_{L^2} \leq \delta $$ the problem \eqref{eq.lwanon} has a
unique global solution
$$ u(t,x) \in C( \R, L^2(\R^n)) \cap Y^\prime.$$
\end{thm}

The proof of Theorem \ref{mainsmfree} is based on the estimate
\eqref{eq.sm} due to Kenig, Ponce, Vega \cite{KPV93},
\cite{KPV98}. In order to have a self contained article we give an
alternative proof of this result due to Kenig, Ponce, Vega in
Section \ref{proofKPV}.

The key step to derive the estimate \eqref{eq.strsmoofree} from
the estimate \eqref{eq.sm} is the following equivalence norm
result.

\begin{thm}\label{comestS.1}
For $n \geq 3, 1 < q < \infty$, for $s\in [-1,1]$ and $a \in \R$
that satisfy
\begin{equation}\label{assS.9}
|a|+|s|<\frac{n}{2},
\end{equation}
the following norms are equivalent
\begin{eqnarray}\nonumber
 \left(\sum_{k \in \mathbb{Z}} 2^{q k a}\|Q_k |D|^{s}f \|^q_{L^2} \right)^{1/q}
 , \\
\label{eqS.10}  \left( \sum_{k \in \mathbb{Z}} 2^{q k a}\|\
|D|^{s} Q_k f \|^q_{L^2} \right)^{1/q} ,
 \\
 \left( \sum_{k \in \mathbb{Z}} \| \
|D|^{s}|x|_k^a f \|^q_{L^2}\right)^{1/q}, \nonumber
\end{eqnarray}
where $|x|^{a}_k = |x|^{a}Q_k(x)$ and the Paley - Littlewood
partition of unity $Q_k(x)$  is defined in \eqref{eq.PLspace}. For
$q = \infty$ the result is still valid with obvious modification
in \eqref{eqS.10}.
\end{thm}

The main idea to establish the  Theorem is similar to the approach
developed in \cite{Gv2000}, \cite{DaGeKu} and \cite{GL2004} for
the case of nonhomogeneous Sobolev spaces and non homogeneous
weights. Therefore, we shall make a localization in coordinate
space and we shall use the Paley Littlewood partition
\eqref{eq.PLspace}. The key point in this approach is to evaluate
the norm of the operator of type $Q_k |D|^{-s} Q_m |D|^s$ with
$|k-m|$ large enough.

The proof of Theorem \ref{mainsmfree}  can be obtained from the
estimate for the Cauchy problem with initial data $f=0$ and the
following Theorems (see section \ref{phaseloc} for the definition
of the spaces $ \ell^{r,s}_D B$ for any Banach space $B$).
\begin{thm}\label{therP.1}
If $q\in[1,2]$ and $a , s \in \R$ satisfy
\begin{equation}\label{tP.7}
\begin{cases}
\  |s|\leq1,  \\
 \  |a|+|s|<\frac{n}{2}
\end{cases}
\end{equation}
then
\begin{equation}\label{tP.8intr}
\| f\|_{\ell_D^{2,0}\ell_x^{q, a}\dot {H}^s}\leq C \|
f\|_{\ell_x^{q, a} \dot {H}^s}.
\end{equation}

\end{thm}

\begin{thm}\label{therP.1du}
If $q\in[2,\infty]$ and $a , s \in \R$ satisfy \eqref{tP.7}, then
\begin{equation}\label{tP.8intrdu}
\| f\|_{\ell_x^{q, a} \dot {H}^s} \leq C \|
f\|_{\ell_D^{2,0}\ell_x^{q, a}\dot {H}^s}.
\end{equation}

\end{thm}

The plan of the work is the following.   The proof
 of the free smoothing estimate of Theorem \ref{mainsmfree} is given in Section \ref{locdef}.
 The proof of the main scale invariant smoothing estimate of
 Theorem \ref{mainsm} is done in Section \ref{profmain}.
 In Section \ref{commest} we treat the commutator estimates needed
 in the
 proof the equivalence of the norms in Theorem \ref{comestS.1}. Some convolution type inequalities needed
 in the proofs of Theorem \ref{comestS.1}
 are included in Section \ref{discr}. The
 concluding steps in the proof of Theorem \ref{comestS.1}
  are presented in Section \ref{SpLoc}.
 Finally the phase localization and the proofs of Theorems
 \ref{therP.1} and \eqref{therP.1du} are given in the last Section \ref{phaseloc}.
 The proof of the estimate
 due to K\"onig, Ponce, Vega is presented in Section \ref{proofKPV} for self contained completeness.

\section{Weighted
 Sobolev spaces estimate of the free Schr\"odinger equation. }

\label{locdef}

Given any Banach space $B\subset D^{\prime}(\R^n)$ satisfying the
property
\begin{eqnarray}\label{S.1a}
\text{for any $Q(x)\in C^\infty_0(\R^n),$} \ \  f\in B \Rightarrow
Q(x)f\in B,
\end{eqnarray}
we can define for any $q\in [1, \infty]$ and for any $\alpha \in
\R$ the space $\ell_x^{q,\alpha}B$ as follows
\begin{equation}\label{S.2a}
\|f\|_{_{\ell_x^{q,\alpha}B}}=\left (\sum_{k\in \Z} \|Q_kf\|_B^q
2^{kq\alpha}\right )^{1/q},
\end{equation}
with obvious modification for $q=\infty$. Note that for any $f\in
C^\infty_0(\R^n \setminus \{0\})$ we have
\begin{equation*}
\|f\|_{\ell_x^{q,\alpha}B}<\infty.
\end{equation*}
So $\ell_x^{q,\alpha}B$ can be defined as the closure of
$C^\infty_0(\R^n \setminus \{0\})$ with respect to the norm
\eqref{S.2a}. An alternative definition is based on the map
\begin{equation}\label{S.3a}
J: f\in C^\infty_0(\R^n\setminus\{0\}) \subset B \rightarrow
J_B(f)_k=\|Q_kf \|_{B}\in \ell^{q,\alpha},
\end{equation}
where $\ell^{q,\alpha}$ is the space of all sequences
$a=(a_k)_{k\in \Z}$ such that
\begin{equation}\label{S.4a}
\|a\|_{\ell^{q,\alpha}}=\left (\sum_{k\in \Z} \|a_k\|^q
2^{kq\alpha}\right )^{1/q} < \infty,
\end{equation}
with obvious modification for $q=\infty$. Then
\begin{equation}\label{S.5a}
\|f\|_{\ell_x^{q,\alpha}B}=\|J_B(f)\|_{\ell^{q,\alpha}}.
\end{equation}
The space $\ell_x^{q,\alpha}B$ is independent of the concrete
choice of Paley-Littlewood decomposition
\begin{equation}\label{S.6a}
\sum_{j\in \Z}Q_j(x)=1
\end{equation}
satisfying
\begin{equation}\label{S:7a}
\begin{cases}
\ Q_j(x)\geq 0,  \\
 \ supp \  Q_j(x) \in \{|x|\sim 2^j\}.
\end{cases}
\end{equation}
A typical example, needed for the smoothing resolvent  type
estimates, is the case $B=\dot{H}^s_p$ where $s\in(-1,1),
 \ 1<p<\infty$. For $s>-\frac{n}{p}$ we have
$\dot{H}^s_p(\R^n)\subset D^\prime (\R^n)$ (see \cite{DaGeKu}) and
the norm is defined by
\begin{equation}\label{S.8a}
\|f \|_{\dot{H}^s_p}=\||D|^sf \|_{L^{p}}.
\end{equation}

After this preparation we can turn to the proof of Theorem
\ref{mainsmfree}.  Starting with the estimate \eqref{eq.sm}, we
use Lemma \ref{LemP.1psin} and find
$$  \sup_{k \in \mathbb{Z}} \||x|_k^{-1/2} \nabla u \|_{L^2_t
L^2_x} \sim  \sup_{k \in \mathbb{Z}} \||x|_k^{-1/2} |D| u
\|_{L^2_t L^2_x}$$ so \eqref{eq.sm} can be rewritten as
\begin{equation}\label{eq.news}
   \sup_{k \in \mathbb{Z}} \||x|_k^{-1/2} |D| u \|_{L^2_t
L^2_x}  \leq  C \left( \sum_{k \in \mathbb{Z}} \||x|_k^{1/2}  F
\|_{L^2_t L^2_x} \right) .
\end{equation}

 Using the fact that $|D|^s$ commutes with
$\Delta$ one can obtain the following consequence of this estimate
\begin{equation}\label{eq.newssig}
   \sup_{k \in \mathbb{Z}} \||x|_k^{-1/2} |D|^{1-\sigma} u \|_{L^2_t
L^2_x}  \leq  C \left( \sum_{k \in \mathbb{Z}} \||x|_k^{1/2}
|D|^\sigma F \|_{L^2_t L^2_x} \right)
\end{equation}
for any $\sigma \in [0,1].$ In particular for $\sigma =1/2$ we get

\begin{equation}\label{eq.newsnn}
   \sup_{k \in \mathbb{Z}} \||x|_k^{-1/2} |D|^{1/2} u \|_{L^2_t
L^2_x}  \leq  C \left( \sum_{k \in \mathbb{Z}} \||x|_k^{1/2}
|D|^{1/2} F \|_{L^2_t L^2_x} \right)
\end{equation}

To this end, we are in position to apply the result of Proposition
\ref{PropS.1} and derive that
\begin{equation}\label{eq.Equi2}
 \sup_{k \in \mathbb{Z}} \||x|_k^{-1/2} |D|^{1/2} u \|_{L^2_t
L^2_x} \sim \|u\|_{L^2_t \ell_x^{\infty, -1/2}\dot{H}^{1/2}_x},
\end{equation}
so
$$
\sup_{k \in \mathbb{Z}} \||x|_k^{-1/2} |D|^{1/2} u \|_{L^2_t
L^2_x} \sim \sup_{k \in \mathbb{Z}} \||x|_k^{-1/2}  u \|_{L^2_t
\dot{H}^{1/2}_x}.$$ In a similar way Proposition \ref{PropS.1}
implies
\begin{equation}\label{eq.Equi1}
\sum_{k \in \mathbb{Z}} \||x|_k^{1/2} |D|^{-1/2} F \|_{L^2_t
L^2_x} \sim \|f\|_{L^2_t \ell_x^{1, 1/2}\dot{H}^{-1/2}_x},
\end{equation}
so
$$
\sum_{k \in \mathbb{Z}} \||x|_k^{1/2} |D|^{-1/2} F \|_{L^2_t
L^2_x} \sim \sum_{k \in \mathbb{Z}} \||x|_k^{1/2}  F \|_{L^2_t
\dot{H}^{-1/2}_x}.$$ The estimate \eqref{eq.newsnn} reads as
\begin{equation}\label{eq.newsnnn}
   \sup_{k \in \mathbb{Z}} \||x|_k^{-1/2}  u \|_{L^2_t
\dot{H}^{1/2}_x}  \leq  C \left( \sum_{k \in \mathbb{Z}}
\||x|_k^{1/2}  F \|_{L^2_t \dot{H}^{-1/2}_x} \right)
\end{equation}
or using the notations of this section (see \eqref{eq.defYYp} and
the definition \eqref{S.5a}) as
\begin{equation}\label{eq.strsmooIa}
    \left\| \int_0^t e^{i(t-s)\Delta} F(s) ds \right\|_{Y^\prime}  \leq  C  \|F \|_{Y} .
\end{equation}
It is easy to derive a similar estimate
\begin{equation}\label{eq.strsmooIadu}
    \left\| \int_t^\infty e^{i(t-s)\Delta} F(s) ds \right\|_{Y^\prime}  \leq  C  \|F \|_{Y}
    ,
\end{equation}
by the aid of \eqref{eq.strsmooIa} and a duality argument for the
quadratic form
$$Q(F,G) = \int\int_{t>s}\langle e^{i(t-s)\Delta}F(s),
G(t)\rangle_{L^2(\R^n)} ds dt.$$ Further, we have to derive the
estimate
\begin{equation}\label{eq.ensmooIa}
    \left\| \int_0^t e^{i(t-s)\Delta} F(s) ds \right\|_{L^\infty_t L^2_x}  \leq  C  \|F \|_{Y} .
\end{equation}
For the purpose set $ u(t) = \int_0^t e^{i(t-s)\Delta} F(s) ds. $
Then $ u=u(t,x)$ is a solution to
\begin{equation}\label{eq.nlfree}
\begin{cases}
\partial_t u - i\Delta u = F,\ \ \  t >0 , \ \ \ x \in \mathbb{R}^n,\\
u(0,x)=0
\end{cases}
\end{equation}
Multiplying by $u$ integration over $\{ 0 \leq t \leq T, x \in
\R^n \infty\}$ we get
$$
\|u(T)\|^2_{L^2(\R^n)} \leq \int_0^T \langle F(t),
u(t)\rangle_{L^2(\R^n)} dt \leq \|F \|_{Y} \|u \|_{Y^\prime} .$$
Applying \eqref{eq.strsmooIa}, we arrive at \eqref{eq.ensmooIa}.
In a similar way we get
\begin{equation}\label{eq.ensmooIadu}
    \left\| \int_0^t e^{i(t-s)\Delta} F(s) ds \right\|_{Y^\prime}  \leq  C  \|F \|_{L^1_t L^2_x} .
\end{equation}
Finally, it remains to prove
\begin{equation}\label{eq.smoohom}
    \left\|  e^{it\Delta} f  \right\|_{Y^\prime}  \leq  C  \|f \|_{ L^2_x} .
\end{equation}
  Consider the
operator $L$ defined by
$$
L : f\in L^{2}_x \  \Longrightarrow \ e^{it\Delta}f. $$ Our goal
is to show that  $L$ is bounded from $L^{2}_x$ to $L^2_t
\ell_x^{\infty, -1/2}\dot{H}^{1/2}_x$. But the continuity of $L$
from $L^{2}_x$ to $L^2_t  \ell_x^{\infty, -1/2}\dot{H}^{1/2}_x$ is
follows from the continuity of its (formally) adjoint
\begin{equation*}
L^*f=\int _0^\infty e^{-i\tau \Delta}f(\tau) d\tau,
\end{equation*}
 from
$Y$ to $L^{2}_x$, which in turns follows from \eqref{eq.ensmooIa}
and the fact that $  e^{it\Delta}$ is unitary operator in $L^2.$

From \eqref{eq.newsnnn}, \eqref{eq.ensmooIa},
\eqref{eq.ensmooIadu}, \eqref{eq.smoohom} and standard energy
estimate, we get \eqref{eq.strsmoomag} and the proof of Theorem
\ref{mainsmfree}  is completed. \\

\section{Proof of Theorem \ref{mainsm} }

\label{profmain}

In this section we will prove the Theorem \ref{mainsm}, so we
shall prove the estimate \eqref{eq.strsmoo}, where $u$ is the
solution of the problem \eqref{eq.lwa}. First of all we have the
identities
\begin{eqnarray}\label{eq.covl2}
    \Delta_A u
   &=& \sum_{j=1}^n (\partial_{x_j} - i A_j) (\partial_{x_j} - i
   A_j)u
   \nonumber\\
   &=&\Delta u -2i\nabla \cdot (Au) + W u ,
\end{eqnarray}
where $$W(t,x) = |A(t,x)|^2  - i \nabla \cdot A$$ satisfies
$$
 \ \sum_{k \in \mathbb{Z}} \ 2^{2k} \ \|W(t,x)\|_{L^\infty_t
L^\infty_{\{|x| \sim 2^k\}}}\leq \varepsilon
$$
due to \eqref{eq.ass2}.

  So,
after a substitution in the equation of \eqref{eq.lwa}, we obtain
\begin{equation}\label{eq.lwa2}
\begin{cases}
i\partial_t u-\Delta u=-2i\nabla\cdot  (Au)+Wu +F  \ \ \  t \in \mathbb{R} , \ \ \ x \in \mathbb{R}^n\\
u(0,x)=f(x).
\end{cases}
\end{equation}
 First of all, we observe that the term $Wu$,
 thanks to the smallness assumption \eqref{ass2}, can be absorbed in the left side of the estimate \eqref{eq.strsmoo}.
 This fact  suggests to localize our attention to the reduced problem
\begin{equation}\label{eq.lwa3}
\begin{cases}
\partial_t u - i\Delta u = -2i\nabla\cdot  (Au) + F,\ \ \  t \in \mathbb{R} , \ \ \ x \in \mathbb{R}^n\\
u(0,x)=f(x).
\end{cases}
\end{equation}
So, using the norms in the spaces $\ell_x^{p,\alpha} B$ introduced
in \eqref{S.2a}   we apply the estimate \eqref{eq.strsmoo}    and
obtain
 \begin{equation}\label{eq.strsm1}
  \|  u  \|_{L^2_t \ell_x^{\infty, -1/2}\dot{H}^{1/2}_x} \leq
  C  \|\nabla\cdot  (Au) \|_{L^2_t \ell_x^{1, 1/2}\dot{H}^{-1/2}_x} + C\| F \|_{L^2_t \ell_x^{1, 1/2}\dot{H}^{-1/2}_x}  + C \|f\|_{L^2}.
 \end{equation}
From  the equivalent norm estimates in Proposition \ref{PropS.1}
(see the equivalence norms relations in \eqref{eq.Equi2},
\eqref{eq.Equi1} also) we have
 \begin{eqnarray}\label{eq.strsm2}
  \| \nabla \cdot (Au) \|_{L^2_t \ell_x^{1, 1/2}\dot{H}^{-1/2}_x}
  \sim\| Au \|_{L^2_t \ell_x^{1, 1/2}\dot{H}^{1/2}_x}.
\end{eqnarray}
From Proposition \ref{PropS.1fg} we have
$$
\| Au \|_{L^2_t \ell_x^{1, 1/2}\dot{H}^{1/2}_x} \lesssim \| A
\|_{L^\infty_t \ell_x^{1, 1}\dot{H}^{1/2}_{2n}} \| u \|_{L^2_t
\ell_x^{\infty, -1/2}L^{2n/(n-1)}_x} + $$ $$ + \| A \|_{L^\infty_t
\ell_x^{1, 1} L^{\infty}_x} \| u \|_{L^2_t \ell_x^{\infty,
-1/2}\dot{H}^{1/2}_x}.
$$
From the Sobolev embedding $ \dot{H}^{1/2}_x \subset
L^{2n/(n-1)}_x, $ we obtain
$$ \| u \|_{L^2_t
\ell_x^{\infty, -1/2}L^{2n/(n-1)}_x} \lesssim \| u \|_{L^2_t
\ell_x^{\infty, -1/2} \dot{H}^{1/2}_x}, $$ while the interpolation
inequality of Proposition \ref{PropS.1fint} guarantees that
$$
\| A \|^2_{L^\infty_t \ell_x^{1, 1}\dot{H}^{1/2}_{2n}} \lesssim \|
\nabla A \|_{L^\infty_t \ell_x^{1, 3/2}L^{2n}} \| A \|_{L^\infty_t
\ell_x^{1, 1/2}L^{2n}}
$$
so applying the H\"older inequality $$\| g \|_{ \ell_x^{1,
a}L^{p}} \lesssim \| g \|_{ \ell_x^{1, a+n/p}L^{\infty}},   $$ we
get
$$
\| A \|^2_{L^\infty_t \ell_x^{1, 1}\dot{H}^{1/2}_{2n}} \lesssim \|
\nabla A \|_{L^\infty_t \ell_x^{1, 2}L^{\infty}} \| A
\|_{L^\infty_t \ell_x^{1, 1}L^{\infty}} \leq \varepsilon^2
$$
due to assumption on $A.$ The above observation implies
$$
\| Au \|_{L^2_t \ell_x^{1, 1/2}\dot{H}^{1/2}_x} \lesssim
\varepsilon \| u \|_{L^2_t \ell_x^{\infty, -1/2}\dot{H}^{1/2}_x}.
$$

Using again the estimates \eqref{eq.strsm1}, we obtain
\begin{equation}\label{eq.strsmoo0I}
     \sup_{k \in \mathbb{Z}} \||x|_k^{-1/2}  u(t,  \cdot) \|_{L^2_t \dot{H}^{1/2}_x}  \leq
  C \| F \|_{L^2_t \ell_x^{1, 1/2}\dot{H}^{-1/2}_x} + C \| f \|_{ L^{2}_x}.
\end{equation}
 This concludes the proof of the theorem.

\section{Application to the semilinear Schr\"odinger equation}

Turning to the semilinear Schr\"odinger equation

\begin{equation}\label{eq.A.1}
\partial_t u - i\Delta_A u = |V u|^{p},
\end{equation}
we note that the class of potentials $V=V(t,x)$, satisfying
\eqref{eq.assV}, obeys certain rescaling property,
 thus one can compute the scaling critical regularity
$$ s = \frac{n}2
- \frac{2-ap}{p-1} $$ and one can expect a well posedness for
initial data $f \in L^2$ if
$$ p = \frac{n+4}{n+2a}.$$
To verify this we shall construct a sequence $u_k(t,x)$ of
functions defined as follows: $u_{-1}(t,x) = 0,$ then we define
the recurrence relation
$$ u_k \rightarrow u_{k+1}(t,x)$$  so that
\begin{equation}\label{eq.lwanonrec}
\begin{cases}
\partial_t u_{k+1} - i\Delta_A u_{k+1} = V(t,x) u_k |u_k|^{p-1},\ \ \  t \in \mathbb{R} , \ \ \ x \in \mathbb{R}^n\\
u_{k+1}(0,x)=f(x).
\end{cases}
\end{equation}
The estimate \eqref{eq.strsmoomag} suggests to show the
convergence of the sequence $u_k$  in the Banach space
$$ Z =  L^\infty_t L^2_x \cap Y^\prime.$$ The definition of the
recurrence relation \eqref{eq.lwanonrec} shows that we have to
show first the property: the map
$$
u \in Z =  L^\infty_t L^2_x \cap Y^\prime \rightarrow V(t,x) u
|u|^{p-1} \in L^1_t L^2_x + Y $$ is a well defined continuous
operator provided $V$ satisfies \eqref{eq.assV}. Our goal is to
show
$$
    \|u_{k+1}\|_{L^\infty_t L^2_x} + \|u_{k+1}\|_{Y^\prime} \leq
    C \|f\|_{L^2} + C \left( \|u_{k}\|_{L^\infty_t L^2_x} +
    \|u_{k}\|_{Y^\prime}\right)^p
$$
or shortly
\begin{equation}\label{eq.A.18}
    \|u_{k+1}\|_{Z}  \leq
    C \|f\|_{L^2} + C  \|u_{k}\|_{Z}^p.
\end{equation}
To apply a contraction argument we need also the inequality
\begin{equation}\label{eq.A.19}
    \|u_{k+1}-u_k\|_{Z}  \leq
     C  \|u_{k}- u_{k+1}\|_{Z}\left(\|u_{k}\|_{Z}+\|u_{k-1}\|_{Z} \right)^{p-1}.
\end{equation}
Combining \eqref{eq.A.18} and \eqref{eq.A.19}, taking
$\|f\|_{L^2}$ sufficiently small, we can show via contraction
argument that $u_k$ converges in $Z$ to the unique solution of
\eqref{eq.A.1} with initial data $u(0)=f.$ Since the proofs of
\eqref{eq.A.18} and \eqref{eq.A.19} are similar, we treat
\eqref{eq.A.18} only. We need actually to verify

\begin{equation}\label{eq.A.20}
\||Vu|^p\|_{Y} \leq C \|u\|^p_{Z}
\end{equation}
To verify this inequality, we start with the definition of the
space $Y$
\begin{equation}\label{eq.A20eq}
    \|g\|_Y  \sim \sum_m  2^{m/2}\| |D|^{-1/2} \varphi \left( \frac{\cdot}{2^m} \right)
    g\|_{L^2_t L^2_x}.
\end{equation}
We apply this relation with $g = |Vu|^p$, combined with the
Sobolev embedding $ \dot{H}^{1/2}(\R^n) \hookrightarrow L^q(\R^n,$
with \begin{equation}\label{eq.A.21} 1/q - 1/2 = 1/(2n),
\end{equation}
and get
\begin{equation}\label{eq.A.22}
\||Vu|^p\|_{Y}  \leq \sum_m  2^{m/2}  \|  \varphi \left(
\frac{\cdot}{2^m} \right)
    Vu\|^p_{L^{2p}_t L^{pq}_x}
\end{equation}
We can apply now the interpolation inequality

\begin{equation}\label{eq.A.23}
\| \phi \|_{L^{2p}_t L^{pq}_x } \leq C \left( \| \phi
\|_{L^{\infty}_t L^{r_1}_x }  \right)^{\theta} \left( \| \phi
\|_{L^{2}_t L^{r_2}_x }  \right)^{1-\theta},
\end{equation}
where $\theta \in (0,1)$ satisfy the relations

\begin{equation}\label{eq.A.24}
\frac{1}{2p} = \frac{1-\theta}{2}, \  \  \frac{1}{pq} =
\frac{\theta}{r_1} + \frac{1-\theta}{r_2}.
\end{equation}
Hence $p(1-\theta) =1.$

From \eqref{eq.A.22} and \eqref{eq.A.23} we get
\begin{equation}\label{eq.A.27}
\||Vu|^p\|_{Y}  \leq \sum_m  2^{m/2}  \|
Vu\|^{p\theta}_{L^{\infty}_t L^{r_1}_{|x|\sim 2^m}} \|
Vu\|^{p(1-\theta)}_{L^{2}_t L^{r_2}_{|x|\sim 2^m}}.
\end{equation}
We choose $r_1 \in (1,2]$ so that  \begin{equation}\label{eq.A.28}
\|V\psi\|_{ L^{r_1}_{|x|\sim 2^m}} \leq C \|\psi\|_{
L^{2}_{|x|\sim 2^m}}
\end{equation}  so taking into account the assumption on
$V$ we see that
\begin{equation}\label{eq.A.29}
   \frac{a}{n} = \frac{1}{r_1} - \frac{1}{2},
\end{equation}
while $r_2 \in (1,2]$ is chosen so that
$$ 2^{m/(2p(1-\theta))}\|V\psi\|_{ L^{r_2}_{|x|\sim 2^m}} \leq C 2^{-m} \|\psi\|_{ L^{2}_{|x|\sim 2^m}} , $$
i.e.
\begin{equation}\label{eq.A.31}
 \frac{1}{r_2} = \frac{1}{2} + \frac{a-1}{n} - \frac{1}{2np(1-\theta)}
 = \frac{1}{2} + \frac{a}{n} - \frac{3}{2n},
\end{equation}
since $p(1-\theta)=1.$  From \eqref{eq.A.24} we find
$$
\frac{1}{pq} = \frac{1}{2} + \frac{a}{n} - \frac{3}{2pn},
$$
so this relation and \eqref{eq.A.21} implies that
$$ p = \frac{n+4}{n+2a}.$$ From \eqref{eq.A.27}, \eqref{eq.A.28}
and \eqref{eq.A.29} we get
$$
\||Vu|^p\|_{Y}  \leq   \| u\|^{p-1}_{L^{\infty}_t L^{2}_{x}} \|
|x|^{-1}u\|_{L^{2}_t L^{2}_x}.
$$
From the Hardy inequality we have
$$ \|
|x|^{-1}u\|_{L^{2}_t L^{2}_x} \leq C \| u\|_{Y^\prime}$$ so
$$
\||Vu|^p\|_{Y}  \leq   \| u\|^{p-1}_{L^{\infty}_t L^{2}_{x}} \|
u\|_{Y^\prime}
$$
and this completes the proof of \eqref{eq.A.20}.

This completes the proof of Theorem \ref{mainsmnonl}.

\section{ESTIMATE OF THE OPERATOR $Q_k |D|^{-s} Q_m |D|^s$.}

\label{commest}

\vspace{0.5cm}

Our goal is to compare  the norms

\begin{equation*}
 \|\ |D|^{-s} f \|_{\ell_x^{q, a} L^p} = \left( \sum_{k
\in \mathbb{Z}} 2^{kq a} \|Q_k |D|^{-s} f\|^q_{L^p} \right)^{1/q}
\end{equation*}

and

\begin{equation*}
 \| f \|_{\ell_x^{q, a} \dot{H}^s_p} = \left( \sum_{k
\in \mathbb{Z}} 2^{kq a} \| |D|^{-s} Q_k f\|^q_{L^p} \right)^{1/q}
\end{equation*}
(see \eqref{S.2a}  or  Section \ref{locdef} for the definition of
the spaces $\ell_x^{q, \alpha}B$, where $B$ is any Banach space
such that $B\subset D^{\prime}(\R^n)$ ).  The key point in the
proof that these norms are equivalent  is the following estimate
of the operator

\begin{equation}\label{E.1}
Q_k |D|^{-s} Q_m |D|^s \ \ \ \text{for} \ \ |k-m|>2 .
\end{equation}

\begin{lem}\label{LemE.1}
For  any $s \in \mathbb{R}, |s|<1,$  any $p, 1<p<n$ and any $k,m
\in \Z , |k-m|\geq 3$ we have the estimate
\begin{equation}\label{E.2}
\|Q_k |D|^{-s} Q_m |D|^s f\|_{L^p} \leq C \  2^{t(k,m,s,p)}\
\|f\|_{L^p},
\end{equation}
where $C=C(s,p)$ independent of $k,m\in \Z$, and

\begin{eqnarray}\label{E.3}
t(k,m,s,p)= k\frac{n}{p}
    + m\frac{n}{p^\prime} - (n-(s\vee0))(k\vee m) -(s\vee
0))(k\wedge m),
\end{eqnarray}
$\frac{1}{p^{\prime}}=1-\frac{1}{p}$,\ $k\wedge m=min(k,m),\ k\vee
m=max(k,m)$.
\end{lem}
{\bf Proof:}
We shall prove the Lemma for $s\in\C$ with $\Re s \in [0,1]$.
For the purpose consider the family of operators
\begin{equation}\label{E.4}
T^{z}=e^{z^2} Q_k |D|^{-z} Q_m |D|^z.
\end{equation}
If $\Re z=0$ then $z=i\sigma,$ $ \sigma \in\mathbb{R}$ and
\begin{equation}\label{E.5}
T^{i\sigma}=e^{-\sigma^2} Q_k |D|^{-i\sigma} Q_m |D|^{i\sigma}.
\end{equation}
Applying stationary phase method (in this case simply integration
by parts), we see that the operator $Q_k |D|^{-i\sigma} Q_m$ has a
kernel
\begin{equation*}
K_{k,m,\sigma}(x,y)
\end{equation*}
satisfying
\begin{equation}\label{E.6}
|K_{k,m,\sigma}(x,y)|\leq C \frac{Q_k(x)Q_m(y)}{2^{n(k\vee m)}}(1+\sigma)^{n+1}.
\end{equation}
This estimate implies
\begin{equation}\label{E.7}
\|Q_k |D|^{-i\sigma}  Q_m  g\|_{L^p} \leq C \frac{2^{k\frac{n}{p}}2^{m\frac{n}{p^{\prime}}}}{2^{n(k\vee m)}}\|g\|_{L^p}.
\end{equation}
Further we apply this inequality with $g=|D|^{i\sigma}f$ and using the following one (see Theorem $1$, Section $2.2$ in \cite{Stei2})
\begin{equation}\label{E.8}
\| |D|^{i\sigma}f\|_{L^{p}}\leq C \| f\|_{L^{p}}(1+\sigma)^{n+1},
\end{equation}
we get
\begin{eqnarray}\label{E.9}
\|T^{i\sigma}(f) \|_{L^{p}}&\leq & C
\frac{2^{k\frac{n}{p}}2^{m\frac{n}{p^{\prime}}}}{2^{n(k\vee m)}}
e^{-\sigma^2} (1+\sigma)^{2(n+1)} \| f\|_{L^{p}}\leq
\nonumber\\
&\leq &C_{1} \frac{2^{k\frac{n}{p}}2^{m\frac{n}{p^{\prime}}}}{2^{n(k\vee m)}} \| f\|_{L^{p}}
\end{eqnarray}
$\forall \sigma\in \R$ with $C_1$ independent of $k,m$ and $\sigma$.
If $z=1+i\sigma$, then
\begin{equation}\label{E.10}
T^{1+i\sigma}=e^{1-\sigma^2+2i\sigma} Q_k |D|^{-1-i\sigma} Q_m \nabla |D|^{i\sigma}\frac{\nabla}{|D|}
\end{equation}
so it is sufficient to estimate the operator
\begin{equation}\label{E.11}
S^\sigma=Q_k |D|^{-1-i\sigma} Q_m \nabla.
\end{equation}
Note that
\begin{equation*}
S^{\sigma}=Q_k (|D|^{-1-i\sigma} \nabla) Q_m-Q_k |D|^{-1-i\sigma} Q^{\prime}_m,
\end{equation*}
where $Q^{\prime}_m=\nabla Q_m$.
The operator $Q_k (|D|^{-1-i\sigma} \nabla) Q_m$ has kernel $K^{\prime}_{k,m}$ satisfying
\begin{equation}\label{E.12}
|K^{\prime}_{k,m}|\leq C\frac{Q_k (x)Q_m(y)}{2^{n(k\vee m)}}
(1+\sigma)^{n+1}
\end{equation}
and this estimate is verified in the same way as \eqref{E.6}. The operator $Q_k |D|^{-1-i\sigma} Q^{\prime}_m$
has kernel  $K^{\prime \prime}_{k,m}$ that satisfies the estimate
\begin{equation}\label{E.13}
|K^{\prime \prime}_{k,m}|\leq C\frac{Q_k (x)Q_m(y)}{2^{(n-1)(k\vee
m)}2^m} (1+\sigma)^{n+1}.
\end{equation}
From \eqref{E.12} and \eqref{E.13} together with \eqref{E.8} we find
\begin{equation}\label{E.14}
\| T^{1+i\sigma}(f)\|_{L^{p}}\leq C \
\frac{2^{k\frac{n}{p}}2^{m\frac{n}{p^{\prime}}}}{2^{(n-1)(k\vee
m)}2^{k\wedge m }}\| f\|_{L^{p}}.
\end{equation}
Applying the complex interpolation argument of Stein (see
\cite{Stei1}), we get \eqref{E.2} for $0<s<1$. This complete the
proof for $s\in(0,1)$.

Next we take $z=-1+i\sigma$. Then
\begin{equation}\label{E.15}
T^{1+i\sigma}=e^{1-\sigma^2-2i\sigma} Q_k |D|^{i\sigma} \frac{\nabla}{|D|} \nabla Q_m |D|^{-1+i\sigma} .
\end{equation}
Then we use the relation
\begin{eqnarray}\label{E.16}
Q_k |D|^{i\sigma} \frac{\nabla}{|D|} \nabla Q_m |D|^{-1+i\sigma}&=&Q_k |D|^{i\sigma} \frac{\nabla}{|D|} Q_m\nabla |D|^{-1+i\sigma}+
\nonumber\\
\\
&&+Q_k |D|^{i\sigma} \frac{\nabla}{|D|} (\nabla Q_m) |D|^{-1+i\sigma}.
\nonumber
\end{eqnarray}
The kernel of $Q_k |D|^{i\sigma} \frac{\nabla}{|D|} (\nabla Q_m)$ is $K^{\prime \prime \prime}_{k,m}$ and satisfies
\begin{equation}\label{E.17}
|K^{\prime \prime \prime}_{k,m}|\leq C\frac{Q_k
(x)Q_m(y)}{2^{(n-1)(k\vee m)}2^m}(1+\sigma)^{n+1},
\end{equation}
then we obtain
\begin{equation}\label{E.18}
\| Q_k |D|^{i\sigma} \frac{\nabla}{|D|} (\nabla Q_m)
|D|^{-1+i\sigma}g\|_{L^{p}}\leq
C\frac{2^{k\frac{n}{p}}2^{m\frac{n}{r^\prime}}}{2^{n(k\vee
m)}2^m}\| g\|_{L^{r}}.
\end{equation}
Taking $g=|D|^{-1}f,$ we get (Hardy-Sobolev)
\begin{equation*}
\||D|^{-1}f \|_{L^{r}}\leq \| f\|_{L^{p}}, \ \ \frac{1}{p}-\frac{1}{r}=\frac{1}{n}.
\end{equation*}
From the fact that $\frac{1}{r}=\frac{1}{p}-\frac{1}{n}$ we have
$1-\frac{1}{r}=\frac{1}{r^\prime}=1-\frac{1}{p}+\frac{1}{n}
=\frac{1}{p^{\prime}}+\frac{1}{n}$ we arrive at
\begin{equation}\label{E.19}
\|Q_k |D|^{i\sigma} \frac{\nabla}{|D|} (\nabla Q_m)
|D|^{-1+i\sigma}f \|_{L^{p}}\leq
C\frac{2^{k\frac{n}{p}}2^{m\frac{n}{p^\prime}}}{2^{n(k\vee m)}}\|f
\|_{L^{p}},
\end{equation}
provided $p<n$. Since
\begin{equation}\label{E.20}
\|Q_k |D|^{i\sigma} \frac{\nabla}{|D|} Q_m\nabla |D|^{-1+i\sigma}
\|_{L^{p}}\leq
C\frac{2^{k\frac{n}{p}}2^{m\frac{n}{p^\prime}}}{2^{n(k\vee m)}}
\|f \|_{L^{p}},
\end{equation}
from \eqref{E.16} and \eqref{E.19} we get
\begin{equation}\label{E.21}
\| T^{1+i\sigma}(f)\|_{L^{p}}\leq
C\frac{2^{k\frac{n}{p}}2^{m\frac{n}{p^\prime}}}{2^{n(k\vee m)}}
\|f \|_{L^{p}}.
\end{equation}
The application of the Stein interpolation argument for $z; \Re
z\in [-1,0]$ combined with the above estimate and \eqref{E.9}
guarantees that  \eqref{E.2} is fulfilled for $s\in (-1,0]$ and
this complete the proof of the Lemma. \hfill$\Box$

It is not difficult to extend the result of Lemma \ref{LemE.1} for
$|k-m|\leq 3$. Note that a formal calculus of $t(k,m,s,p)$ for
$|k-m|\leq 3$ in \eqref{E.3} gives $2^{t(k,m,s,p)}\sim 1$. To
verify
\begin{equation}\label{E.22}
\|Q_k |D|^{-s}  Q_m |D|^{s}f \|_{L^{p}}\leq C \|f \|_{L^{p}},
\end{equation}
for $|s|<1, 1<p<n$, it is sufficient to use a scale argument and to show \eqref{E.22} for $k=m=0$ so we shall verify the inequality
\begin{equation}\label{E.23}
\|Q_0 |D|^{-s}Q_0 |D|^{s}f \|_{L^{p}}\leq C\|f \|_{L^{p}}.
\end{equation}
Here we can use an interpolation argument as in the proof of Lemma
\ref{LemE.1}. Then we have to show that
$L_\sigma=Q_0|D|^{-1+i\sigma}Q_0\nabla$ is $L^{p}$-bounded. But
\begin{equation}\label{E.24}
L_\sigma=Q_0(|D|^{-1+i\sigma}\nabla )Q_0+Q_0|D|^{i\sigma}|D|^{-1}(\nabla Q_0).
\end{equation}
Since $|D|^{-1}\nabla$ is operator of order $0$ it is
$L^p$-bounded and $|D|^{i\sigma}$ is also $L^p$ -bounded, we see
that $Q_0(|D|^{-1+i\sigma}\nabla )Q_0$ is $L^p$ -bounded. From he
property
\begin{equation}\label{E.25}
|D|^{-1}: L^r\rightarrow L^p, \ \  \frac{1}{r}=\frac{1}{p}-\frac{1}{n},
\end{equation}
and
\begin{equation}\label{E.26}
\nabla Q_0 : L^p\rightarrow L^r,
\end{equation}
we see that $ |D|^{-1}(\nabla Q_0): L^p\rightarrow L^p$ so
$L_\sigma$ is $L^p$ -bounded. This observation and a Stein
interpolation argument implies $\eqref{E.22}$ for $s\in (0,1)$. To
cover the case $s\in (-1,0)$ we have to show that
\begin{equation}\label{E.27}
L^{\prime}_\sigma=Q_0|D|^{i\sigma}\frac{\nabla }{|D|}\nabla Q_0 |D|^{-1+i\sigma}, \ \ \  \text{is $L^p$-bounded}
\end{equation}
But
\begin{equation}\label{E.28}
L^{\prime}_\sigma=Q_0 |D|^{i\sigma} \frac{\nabla }{|D|}(\nabla Q_0) |D|^{-1} |D|^{i\sigma}+Q_0|D|^{i\sigma} \frac{\nabla }{|D|} \nabla Q_0 \frac{\nabla}{|D|} |D|^{i\sigma},
\end{equation}
and again we can show that  \eqref{E.25} and \eqref{E.26} imply
that the operator in the right side of the \eqref{E.28} is
$L^p$-bounded ( with norm $\leq C(1+\sigma)^{n+1}$ ). Since the
second operator is also $L^p$ -bounded, we see that
$L^{\prime}_\sigma$ is also $L^p$ -bounded and this argument
implies

\begin{lem}\label{LemE.2}
For any $s\in \R,$  $|s|<1$  any $p, 1<p<n$ there exists a
constant $C=C(s,p,n)>0$ so that for any $k,m\in \Z$, and for $f\in
S(\R^n)$ we have
\begin{equation}\label{E.29}
\|Q_k |D|^{-s} Q_m |D|^s f\|_{L^p} \leq C \  2^{t(k,m,s,p)}\
\|f\|_{L^p},
\end{equation}
 where $t(k,m,s,p)$ is defined  in \eqref{E.3}.
 \end{lem}
Finally we use a duality argument and find :

\begin{lem}\label{LemE.3}
For any $s\in \R$ such that $|s|<1$ we have for any $p, (\frac{n}{n-1}<p<n)$ there exists a constant $C=C(s,p,n)>0$
so that for $f\in S(\R^n)$ we have
\begin{equation}\label{E.30}
\|Q_k |D|^{-s} Q_m |D|^s f\|_{L^p} \leq C \  2^{t(k,m,s,p)}\
\|f\|_{L^p},
\end{equation}
 where $t(k,m,s,p)$ is defined  in \eqref{E.3}.
 \end{lem}
{\bf Proof:} For any$f,g\in S(\R^n)$ we have
\begin{eqnarray}\label{E.31}
|(g,|D|^{-s} Q_k |D|^{s}Q_m f)| &=&|(Q_m |D|^{-s}Q_k|D|^{s}g, f)|
\nonumber\\
&\leq&\|f\|_{L^p}\||D|^{-s} Q_k |D|^{s}Q_m g\|_{L^{p^{\prime}}}
\end{eqnarray}
Applying for $p^{\prime}$  the estimate of \eqref{LemE.2}, we find
\begin{equation}\label{E.32}
\||D|^{-s} Q_k |D|^{s}Q_m g\|_{L^{p^{\prime}}}\leq
2^{t(m,k,s,p^\prime)}\ \|g\|_{L^p},
\end{equation}
where
$$t(m,k,s,p^\prime)=k\frac{n}{p}+ m\frac{n}{p^{\prime}} -(n-(s\vee0))(k\vee
m)- (s\vee 0))(k\wedge m)=t(k,m,s,p).$$ This complete the proof.
\hfill$\Box$

\section{Discrete Estimates.}

\label{discr}

Consider the operator
\begin{equation}\label{D.1}
T:a=\{a_k\}_{k\in \Z}\rightarrow Ta=b_{m}=\sum_{k\in\{|k-m|\geq4\}} t_{k,m}a_k,
\end{equation}
where
\begin{equation}\label{D.2}
t_{k,m}=2^{m\lambda}\ 2^{\mu k}\ 2^{-\beta(m\vee k)}, \ \ m\vee
k=max(m,k),
\end{equation}
\begin{equation}\label{D.3}
\lambda >0, \mu>0,\ \  \beta=\lambda+\mu.
\end{equation}

\begin{lem}\label{LemD.1}
If $\lambda, \mu>0,$ and $\beta=\lambda+\mu,$ then the operator
\begin{equation*}
T: \ell^q\rightarrow \ell^q
\end{equation*}
is bounded for any $q\in[1,\infty]$.
\end{lem}
{\bf Proof:} First we consider the cases $q=\infty$ and $q=1$, then we apply the interpolation argument.
We represent $Ta$ as
\begin{equation}\label{D.4}
Ta=T_1a+T_2a,
\end{equation}
where
\begin{equation}\label{D.5}
(T_1a)_m=\sum_{k=m+1}^{\infty} t_{k,m}a_k,
\end{equation}
\begin{equation}\label{D.6}
(T_2a)_m=\sum_{-\infty}^{m} t_{k,m}a_k.
\end{equation}
From \eqref{D.2} we find for $T_1a$
\begin{eqnarray}\label{D.7}
\ \| T_1a\|_{\infty}\leq Csup_{m\in \Z}\left(
\sum^{\infty}_{k=m+1}t_{k,m}\right) \|a\|_\infty \leq
\nonumber\\
\\
\ \leq C sup_{m\in \Z} \left (
\sum^{\infty}_{k=m+1}2^{m\lambda}2^{k\mu}2^{-\beta m}\right)
\|a\|_\infty. \nonumber
\end{eqnarray}
so
\begin{equation}\label{D.8}
\ \| T_1a\|_{\infty}\leq C \|a\|_\infty.
\end{equation}
From \eqref{D.2}, for $T_2a$ we have the following estimate
\begin{eqnarray}\label{D.9}
\ \| T_2a\|_{\infty}\leq C sup_{m\in \Z} \left (
\sum^{m}_{k=-\infty}2^{m\lambda}2^{\mu k}2^{-\beta k}\right)
\|a\|_\infty =
\nonumber\\
\\ =C sup_{m\in \Z}\left (  \sum^{m}_{k=-\infty}2^{m\lambda}2^{-\lambda k}\right) \|a\|_\infty =C\|a\|_{\infty},
\nonumber
\end{eqnarray}
This estimate and \eqref{D.8} imply
\begin{equation}\label{D.11}
\ \| Ta\|_{\infty}\leq C \|a\|_\infty.
\end{equation}
For $q=1$ we have
\begin{eqnarray}\label{D.12}
\ \| T_1a\|_{1}\leq C sup_{k\in \Z}\left( \sum_{\{m\in\Z;m\geq k\}}t_{k,m}\right)  \|a\|_1\leq
\nonumber\\
\leq  sup_{k\in \Z}\left (  \sum^{\infty}_{m=k}2^{m\lambda}2^{k\mu}2^{-\beta m}\right) \|a\|_1=\\
\ =Csup_{k\in \Z}2^{k\mu} \left (  \sum^{\infty}_{m=k}2^{-m\mu}\right) \|a\|_1\leq 2C \|a\|_1.
\nonumber
\end{eqnarray}
In a similar way we estimate $T_2a$,
\begin{eqnarray}\label{D.13}
\ \ \| T_2a\|_{1}\leq C sup_{k\in \Z}\left (  \sum^{k-1}_{m=-\infty}2^{m\lambda}2^{k\mu}2^{-\beta k}\right) \|a\|_1=
\nonumber\\
\\
 =C sup_{k\in \Z}\left (  \sum^{k-1}_{m=-\infty}\right) 2^{m\lambda}2^{-k\lambda}\|a\|_1\leq C \|a\|_1.
\nonumber
\end{eqnarray}
Thus we get
\begin{equation}\label{D.14}
\ \| Ta\|_{1}\leq C \|a\|_1,
\end{equation}
and this completes the proof of the Lemma.
\hfill$\Box$\\
It easy to obtain the corresponding weighted version of Lemma
\ref{LemD.1} in terms of weighted $\ell^q$ spaces
\begin{equation}\label{D.15}
\ell^{q,\alpha}=\{a=(a)_{k\in \Z}; \sum_k
2^{kq\alpha}|a_k|^q<\infty \}.
\end{equation}
For the purpose consider the operator
\begin{equation*}
J^\alpha: a\rightarrow b=J^\alpha a,
\end{equation*}
defined as follows
\begin{equation}\label{D.16}
b_k=2^{k\alpha} a_k,
\end{equation}
we have the two following Lemmas
\begin{lem}\label{LemD.2}
The application $J^\alpha: \ell^q \rightarrow \ell^{q,\alpha}$ is
an isomorphism for any $\alpha\in\R$ and any $q\in [1, \infty]$.
\end{lem}
\begin{lem}\label{LemD.3}
If $\sigma, \nu, \lambda, \mu$ are real numbers such that
\begin{equation}\label{D.17}
\begin{cases}
\  \lambda+\sigma>0,\\
\ \mu- \nu>0,
\end{cases}
\end{equation}
then for $\beta=\lambda+\sigma+\mu- \nu$ we have
\begin{equation*}
T: \ell^{q,\sigma}\rightarrow \ell^{q,\nu}
\end{equation*}
where $T$ is defined by \eqref{D.1} and \eqref{D.2}.
\end{lem}
{\bf Proof:} Let
\begin{equation*}
\widetilde{T}: J_{\sigma}TJ_{\nu}^{-1}.
\end{equation*}
Then Lemma \ref{LemD.2} guarantees that $T:
l^{q,\sigma}\rightarrow \ell^{q,\nu}$ if and only if
$\widetilde{T}: \ell^{q}\rightarrow \ell^{q}$. Note that
$\widetilde{T}$ by
\begin{equation}\label{D.18}
t_{m,k}=2^{m(\lambda+\sigma)}2^{k(\mu-\nu)}2^{-\beta(m\vee k)}.
\end{equation}
So applying Lemma \ref{LemD.1} with
$\widetilde{\lambda}=\lambda+\sigma$ and
$\widetilde{\mu}=\mu-\nu$, we complete the proof.

A slight generalization of Lemma \ref{LemD.1} can be obtained for
the case when $\lambda, \mu, \beta$ are vectors in $\R^2$, that is
\begin{equation*}
\begin{cases}
\ \lambda=(\lambda_1, \lambda_2) \\
\ \mu=(\mu_1, \mu_2)\\
\ \beta=(\beta_1, \beta_2).
\end{cases}
\end{equation*}
Then \eqref{D.1}
\begin{equation}\label{D.19}
\begin{cases}
\ Ta=b, \ \ \text{where} \ \ b_m=\sum_{k\in \Z^2} t_{mk}a_k, \  m\in\Z^2\\
\ a=(a)_{k\in \Z^2},
\end{cases}
\end{equation}
where
\begin{equation}\label{D.20}
t_{m,k}=2^{{\sum_{j=1}^{2}m_j \lambda_j+k_j\mu_j-\beta_j(m_j\vee
k_j)}}.
\end{equation}
The assumption \eqref{D.2} can be replaced again by the following one
\begin{equation}\label{D.21}
\lambda_j>0, \mu_j>0, \ j=1,2.
\end{equation}
\begin{lem}\label{LemD.4}
If $\lambda, \mu, \beta\in \R^2$
satisfy $\beta_j=\lambda_j+\mu_j, i=1,2$ and \eqref{D.18} then
\begin{equation}\label{D.22}
T: \ell^{q_1}_{k_1} \ell^{q_2}_{k_2}\rightarrow \ell^{q_1}_{k_1}
\ell^{q_2}_{k_2},
\end{equation}
is bounded for $q=(q_1,q_2), 1\leq q_j\leq \infty$.
\end{lem}
\begin{remark}\label{D.23}
Given any sequence $a=\{a_{k_1k_2}\}_{k={k_1,k_2}\in \Z}$
we can consider the norm
\begin{equation}\label{D.24}
\| a\|_{\ell^{q_1}_{k_1} \ell^{q_2}_{k_2}}=\left (\sum_{k_1\in \Z}
\left (\sum_{k_2\in \Z}| a_{k_1k_2}|^{q_2}\right )^{q_1/q_2}\right
)^{1/q_{1}},
\end{equation}
( with obvious modification if $q_1=\infty$ or $q_2=\infty$ ), and
the corresponding Banach space $\ell^{q_1}_{k_1}
\ell^{q_2}_{k_2}$. Note that $$\ell^{q_1}_{k_1} \ell^{q_2}_{k_2}
\ne \ell^{q_2}_{k_2}\ell^{q_1}_{k_1},$$ but the assertion of Lemma
\ref{LemD.4} is still true if we replace $\ell^{q_1}_{k_1}
\ell^{q_2}_{k_2}$ by $\ell^{q_2}_{k_2}\ell^{q_1}_{k_1}$. The
corresponding generalization of Lemma \ref{LemD.3} is the
following,
\begin{lem}\label{LemD.5}
If $\sigma, \nu, \lambda, \mu \in \R^2$ satisfy
\begin{equation}\label{D.25}
\begin{cases}
\  \lambda_j+\sigma_j>0,\\
\ \mu_j- \nu_j>0,
\end{cases}
\end{equation}
then for $\beta_j=\lambda_j+\sigma_j+\mu_j- \nu_j, j=1,2$ the
operator  $T$ defined by \eqref{D.16} and \eqref{D.17} is in
$B(\ell^{q_1, \sigma_1}_{k_1} \ell^{q_2,\sigma_2}_{k_2},\ell^{q_1,
\nu_1}_{k_1} \ell^{q_2, \nu_2}_{k_2})$.
\end{lem}
\end{remark}

\section{Space localization.}\label{SpLoc}

Given any Banach space $B\subset D^{\prime}(\R^n)$ satisfying the property
\begin{eqnarray}\label{S.1}
\text{for any $Q(x)\in C^\infty_0(\R^n),$} \ \  f\in B \Rightarrow Q(x)f\in B,
\end{eqnarray}
we can define for any $p\in [1, \infty]$ and for any $a \in \R$
the space $\ell_x^{q,a}B$ as follows
\begin{equation}\label{S.2}
\|f\|_{\ell_x^{q,a} B}=\left (\sum_{k\in \Z} 2^{kqa} \
\|Q_kf\|_B^q \right )^{1/q},
\end{equation}
with obvious modification for $q=\infty$. Note that for any $f\in
C^\infty_0(\R^n\setminus \{0\})$ we have
\begin{equation*}
\|f\|_{\ell_x^{q,a} B}<\infty.
\end{equation*}
So $\ell_x^{q,a} B$ can be defined as the closure of
$C^\infty_0(\R^n\setminus \{0\})$ with respect to the norm
\eqref{S.2}. An alternative definition is based on the map
\begin{equation}\label{S.3}
J: f\in C^\infty_0(\R^n \setminus \{0\}) \subset B \rightarrow
J_B(f)_k=\|Q_kf \|_{B}\in \ell^{q,\alpha},
\end{equation}
where $\ell^{q,\alpha}$ is the space of all sequences
$a=(a_k)_{k\in \Z}$ such that
\begin{equation}\label{S.4}
\|a\|_{l^{q,\alpha}}=\left (\sum_{k\in \Z} \|a_k\|^q 2^{kq\alpha}\right )^{1/q},
\end{equation}
with obvious modification for $p=\infty$. Then
\begin{equation}\label{S.5}
\|f\|_{l^{q,\alpha}B}=\|J_B(f)\|_{l^{q,\alpha}}.
\end{equation}
The space $l^{q,\alpha}B$ is independent of the concrete choice of
Paley-Littlewood decomposition
\begin{equation}\label{S.6}
\sum_{j\in \Z}Q_j(x)=1
\end{equation}
satisfying
\begin{equation}\label{S:7}
\begin{cases}
\ Q_j(x)\geq 0,  \\
 \ supp  \ Q_j(x) \in \{|x|\sim 2^j\}.
\end{cases}
\end{equation}
A typical example is the case $B=\dot{H}^s_p,$ where $s \in(-1,1),
1<p<\infty$. For $s>-\frac{n}{p}$ we have
$\dot{H}^s_p(\R^n)\subset D^\prime$ (see \cite{DaGeKu}) and the
norm is defined by
\begin{equation}\label{S.8}
\|f \|_{\dot{H}^s_p}=\||D|^sf \|_{L^{p}}.
\end{equation}

Our next goal is to show the equivalence of the norms
\begin{equation*}
 \|\ |D|^{-s} f \|_{\ell_x^{q, a} L^p} = \left( \sum_{k
\in \mathbb{Z}} 2^{kq a} \|Q_k |D|^{-s} f\|^q_{L^p} \right)^{1/q}
\end{equation*}
and
\begin{equation*}
 \| f \|_{\ell_x^{q, a} \dot{H}^s_p} = \left( \sum_{k
\in \mathbb{Z}} 2^{kq a} \| |D|^{-s} Q_k f\|^q_{L^p} \right)^{1/q}
\end{equation*}
The proof of the equivalence norm Theorem \ref{comestS.1}  is a
direct consequence (taking $p=2$) of the following estimates.
\begin{prop}\label{PropS.1}
For $p\in (n/(n-1), n), q\in [1, \infty]$, for $s\in [-1,1]$ and
$a\in \R$ that satisfy
\begin{equation}\label{S.9}
|a|+|s|<min \left (\frac{n}{p}, \frac{n}{p\prime} \right)
\end{equation}
one can find a constant $C=C(n,s,p,q,a)>0$ so that
\begin{eqnarray}\label{S.10}
 C^{-1} \|\ |D|^{-s} f \|_{\ell_x^{q, a} L^p} \leq  \| f \|_{\ell_x^{q, a} \dot{H}^s_p} \leq
 C\|\ |D|^{-s} f \|_{\ell_x^{q, a} L^p}.
\end{eqnarray}
\end{prop}
{\bf Proof:} The left inequality in \eqref{S.10} is equivalent to
\begin{equation}\label{S.11}
 \left(\sum_{k \in \mathbb{Z}} 2^{kq a} \|Q_k |D|^{-s} f\|^q_{L^p}
\right)^{1/q} \leq C  \left( \sum_{k \in \mathbb{Z}} 2^{kq a} \|
|D|^{-s} Q_k f\|^q_{L^p} \right)^{1/q}.
\end{equation}
Indeed, given any integers $k,m\in \Z$ with $|k-m|>2$ we have the identity
\begin{equation}\label{S.12}
 Q_k |D|^{-s}Q_mf=Q_k |D|^{-s}Q_m |D|^{s}|D|^{-s}\widetilde{Q}_m f,
\end{equation}
where $\widetilde{Q}_m=\frac{1}{3}(Q_{m-1}+Q_m+Q_{m+1})$ is
another Paley-Littlewood partition of unity $$\sum_{m\in
\Z}\widetilde{Q}_m=1,$$ such that $\widetilde{Q}_m(s)=1$ for $s\in
supp \ Q_m$. To verify \eqref{S.11} it is sufficient to show that
\begin{equation}\label{S.13}
 \left(\sum_{k \in \mathbb{Z}} 2^{kq a} \|Q_k |D|^{-s} f\|^q_{L^p}
\right)^{1/q} \leq C  \left( \sum_{m \in \mathbb{Z}} 2^{mq a} \|
|D|^{-s} \widetilde{Q}_m f\|^q_{L^p} \right)^{1/q}.
\end{equation}
From the estimate of Lemma \ref{LemE.1} we have
\begin{equation}\label{S.14}
\|Q_k |D|^{-s} Q_m |D|^s f\|_{L^p} \leq C 2^{t(k,m,s,p)}\
\|f\|_{L^p},
\end{equation}
where $t(k,m,s,p)$ is defined in \ref{E.3}. Applying the above
estimate with $$g=|D|^{-s}\widetilde{Q}_m f$$ together with  Lemma
\ref{LemD.3}, we complete the proof of \eqref{S.13}.

To verify the right inequality in \eqref{S.11} it sufficient to
show
\begin{equation}\label{S.16}
 \left( \sum_{k \in \mathbb{Z}} 2^{kq a} \|
|D|^{-s} Q_k f\|^q_{L^p} \right)^{1/q} \leq C \left(\sum_{k \in
\mathbb{Z}} 2^{kq a} \|Q_k |D|^{-s} f\|^q_{L^p} \right)^{1/q}.
\end{equation}
To this end we use the relation
\begin{eqnarray*}
|D|^{-s}Q_kf&=&|D|^{-s}Q_k  |D|^{s}|D|^{-s}f=\\
&=&\sum_{m\in \Z} |D|^{-s}Q_k |D|^{s}Q_m \widetilde{Q}_m|D|^{-s}.
\end{eqnarray*}
From Lemma \ref{LemE.3} we have
\begin{equation}\label{S.17}
\|Q_k |D|^{-s} Q_m |D|^s f\|_{L^p} \leq  C 2^{t(k,m,s,p)} \
\|f\|_{L^p} ,
\end{equation}
 where $t(k,m,s,p)$ is defined in \ref{E.3}, so applying Lemma \ref{LemD.3}, we obtain \eqref{S.16} and complete
 the proof of the Proposition.

\hfill$\Box$

\begin{prop}\label{PropS.1fg}
For $p\in (n/(n-1), n), q\in [1, \infty]$, for $s\in [0,1]$ and
$a>0$ that satisfy \eqref{S.9} one can find a constant
$C=C(n,s,p,q,a)>0$ so that
\begin{equation}\label{S.10fg}
  \| f g\|_{\ell_x^{q, a} \dot{H}^s_p} \leq
 C \left(\| f\|_{\ell_x^{q_1, a_1} \dot{H}^{s}_{p_1}} \| g\|_{\ell_x^{q_2, a_2} L^{p_2}}
 + \| f\|_{\ell_x^{q_3, a_3} L^{p_3}} \| g\|_{\ell_x^{q_4, a_4} \dot{H}^{s}_{p_4}
 }\right),
\end{equation}
provided $ a_1,a_2,a_3,a_4 \geq 0$ and $1 \leq p_1,p_2,p_3,p_4,
q_1,q_2,q_3,q_4 \leq \infty$ satisfy
$$ a_1+a_2=a_3+a_4=a, \frac{1}{p} = \frac{1}{p_1} + \frac{1}{p_2} = \frac{1}{p_3} + \frac{1}{p_4},
\frac{1}{q} = \frac{1}{q_1} + \frac{1}{q_2} = \frac{1}{q_3} +
\frac{1}{q_4}.
$$
\end{prop}
{\bf Proof:} The proof uses the previous Proposition and the
standard multiplicative Sobolev inequality
\begin{equation}\nonumber
  \| f g\|_{\dot{H}^s_p} \leq
 C \left(\| f\|_{ \dot{H}^{s}_{p_1}} \| g\|_{ L^{p_2}}
 + \| f\|_{ L^{p_3}} \| g\|_{ \dot{H}^{s}_{p_4}
 }\right),
\end{equation}
so we omit the details.

\hfill$\Box$

Using the interpolation property
$$
(H^{s_1}_{p_1}, H^{s_1}_{p_1})_\theta =  H^s_p, \  s=(1-\theta)s_1
+ \theta s_2, \ \frac{1}{p} = \frac{1-\theta}{p_1} +
\frac{\theta}{p_2},
$$
we arrive at

\begin{prop}\label{PropS.1fint}
For $p\in (n/(n-1), n), q\in [1, \infty]$, for $s\in (0,1)$ and
$a>0$ that satisfy \eqref{S.9} one can find a constant
$C=C(n,s,p,q,a)>0$ so that
\begin{equation}\label{S.10fgint}
  \| f\|_{\ell_x^{q, a} \dot{H}^s_p} \leq
 C \left(\| f\|_{\ell_x^{q_1, a_1} \dot{H}^{1}_{p_1}}\right)^{1-\theta} \left( \| f\|_{\ell_x^{q_2, a_2} L^{p_2}}
  \right)^\theta,
\end{equation}
provided $ a_1,a_2 \geq 0$ and $1 < p_1,p_2,q_1,q_2 <\infty$
satisfy
$$a = a_1(1-\theta)+a_2\theta, s = 1-\theta, \frac{1}{p} = \frac{1-\theta}{p_1} + \frac{\theta}{p_2} ,
\frac{1}{q} = \frac{1-\theta}{q_1} + \frac{\theta}{q_2}.
$$
\end{prop}

\begin{remark}\label{MainEq}
Note that, using the norms introduced in \eqref{S.2}, we the
estimate \eqref{eq.strsmoofree} can be written as
\begin{equation}\label{eq.strsmoofreem}
 \|  u  \|_{L^2_t \ell_x^{\infty, -1/2}\dot{H}^{1/2}_x} \leq C \| f \|_{ L^{2}_x} + C  \| F  \|_{L^2_t \ell_x^{1, 1/2}\dot{H}^{-1/2}_x} .
\end{equation}

\end{remark}

\section{Phase localization}

\label{phaseloc}

Given any Banach space $B\subset D^{\prime}(\R^n)$ satisfying the
property
\begin{eqnarray}\label{P.1}
\text{for any $P(\xi)\in C^\infty_0(\R^n),$} \ \  f\in B
\Rightarrow P(D)f\in B,
\end{eqnarray}
we can define for any $r\in [1, \infty]$ and for any $s\in \R$ the
space $\ell^{r,s}_D B$ as follows
\begin{equation}\label{P.2}
\|f\|_{\ell_D^{r,s}B}=\left (\sum_{k\in \Z} \|P_k(D) f\|_B^r
2^{krs}\right )^{1/r},
\end{equation}
with obvious modification for $r=\infty$. Here  $\{P_k(\xi)\}$ is
a Paley-Littlewood decomposition.

 Our goal is to find some
concrete examples of Banach spaces $B$ satisfying the embedding
\begin{equation}\label{P.3}
B\subset \ell_D^{r,0}B.
\end{equation}

Therefore  we look for estimate of type
\begin{equation}\label{P.4}
\left (\sum_{k\in \Z} \| P_k(D)f\|_B^r \right )^{1/r}\leq C\|
f\|_{B}.
\end{equation}

A typical example for Banach space $B$ satisfying \eqref{P.4} is
$B=L^p$ with $1<p\leq 2$, so
\begin{equation}\label{P.4a}
 \left (\sum_{k\in \Z} \|P_k(D)f \|_{L^p} ^2 \right )^{1/2}\leq C\| f\|_{L^{p}}, \ \ 1<p\leq 2.
\end{equation}
Having in mind that spaces $\ell_x^{r,a}\dot {H}^s_p$ are natural
candidate for estimate of type \eqref{P.4}, we shall verify that
the conditions
\begin{equation}\label{P.6}
\begin{cases}
\ 1\leq q\leq 2, \ p=2  \\
 \  |a|+|s|<\frac{n}{2},
\end{cases}
\end{equation}
imply \eqref{P.4}. More precisely we have
\begin{lem}\label{LemP.1}
If $q\in[1,2]$ and $a , s \in \R$ satisfy
\begin{equation}\label{P.7}
\begin{cases}
\  |s|\leq1,  \\
 \  |a|+|s|<\frac{n}{2},
\end{cases}
\end{equation}
then
\begin{equation}\label{P.8}
\| f\|_{\ell_D^{2,0}\ell_x^{q,a}\dot {H}^s}\leq C \|
f\|_{\ell_x^{q,a} \dot {H}^s}.
\end{equation}

\end{lem}

{\bf Proof:} For any $f\in S(\R^n)$ we have ( for any $r,q\in (1, \infty)$)
\begin{eqnarray}\label{P.9}
 \| f\|_{\ell_D^{2,0}\ell_x^{q,a}\dot {H}^s}\cong
 \left(\sum_{k_2 \in \mathbb{Z}} \left( \sum_{k_1 \in \mathbb{Z}} 2^{k_1 q a}\||D|^s Q_{k_1} P_{k_{2}}(D)f \|^q_{L^2}\right)^{2/q} \right)^{1/2}   \cong
 \\
\ \
 \cong \left(\sum_{k_2 \in \mathbb{Z}}
 \left( \sum_{k_1 \in \mathbb{Z}} 2^{k_1 q a}\| Q_{k_1} |D|^s P_{k_{2}}(D)f \|^q_{L^2}\right)^{2/q} \right)^{1/2}
 \cong \| \|Q_{k_1} P_{k_{2}}(D)f \|_{L^2}\|_{\ell_{k_2}^{2,s} \ell_{k_1}^{q,a}}, \nonumber
\end{eqnarray}
where here and below we use the discrete norm in $ l_{k_1}^{q,a}
 l_{k_2}^{2,s}$  introduced in
\eqref{D.24}. In the second equivalence relation we have used
Proposition \ref{PropS.1}. Further we have
\begin{equation}\label{P.10}
Q_{k_1} P_{k_{2}}(D)f=\sum_{m_{1}\in \Z} \sum_{m_{2}\in \Z}
Q_{k_1} P_{k_{2}}(D)Q_{m_1} P_{m_{2}}(D)\widetilde
{P}_{m_{2}}(D)\widetilde {Q}_{m_1} f.
\end{equation}

It is not difficult, using again integration by parts argument, to see that for $|k_1-m_1|\geq 3$ we have
\begin{equation}\label{P.11}
\| Q_{k_1} P_{k_{2}}(D)Q_{m_1}g\|_{L^{2}}\leq C \frac{2^{(k_1+m_1)n/2}}{2^{(k_1\vee m_1)n}}\| g\|_{L^{2}}
\end{equation}
so
\begin{equation}\label{P.12}
\| Q_{k_1} P_{k_{2}}(D)Q_{m_1}P_{m_{2}}(D)f\|_{L^{2}}\leq C \frac{2^{(k_1+m_1)n/2}}{2^{(k_1\vee m_1)n}}\| f\|_{L^{2}}.
\end{equation}
In a similar way, using the same integration by parts argument, we find for $|k_2-m_2|\geq 3$
\begin{equation}\label{P.13}
\|  P_{k_{2}}(D)Q_{m_1}P_{m_{2}}(D)g\|_{L^{2}}\leq C \frac{2^{(k_2+m_2)n/2}}{2^{(k_2\vee m_2)n}}\| g\|_{L^{2}},
\end{equation}
so
\begin{equation}\label{P.14}
\| Q_{k_1} P_{k_{2}}(D)Q_{m_1}P_{m_{2}}(D)f\|_{L^{2}}\leq C \frac{2^{(k_2+m_2)n/2}}{2^{(k_2\vee m_2)n}}\| f\|_{L^{2}}.
\end{equation}
An interpolation between \eqref{P.12} and  \eqref{P.14} gives
\begin{equation}\label{P.15}
\| Q_{k_1} P_{k_{2}}(D)Q_{m_1}P_{m_{2}}(D)f\|_{L^{2}}\leq C t_{k,m}^{(\theta)}\| f\|_{L^{2}}.
\end{equation}
where $k=(k_1,k_2)\in \Z^2, m=(m_1,m_2)\in \Z^2$, $\theta \in [0,1] $ will be chosen later on and
\begin{equation}\label{P.16}
 t_{k,m}^{(\theta)}= \frac{2^{(k_1+m_1)\theta n/2}}{2^{(k_1\vee m_1)\theta n}}
 \frac{2^{(k_2+m_2)(1-\theta) n/2}}{2^{(k_2\vee m_2)(1-\theta) n}}.
\end{equation}
If $a, s\in\R$ satisfy

\begin{equation}\label{P.17}
|a|+|s|<\frac{n}{2},
\end{equation}
then we can choose $\theta \in [0,1] $ so that
\begin{equation}\label{P.18}
\begin{cases}
\ |a| \leq \frac{n}{2}(1-\theta),  \\
 \  |s|<\frac{n}{2}\theta .
\end{cases}
\end{equation}
Using the argument of the proof of Lemma \ref{E.2}, we see that
\eqref{P.14} is fulfilled without the restrictions $|k_1-m_1|\geq
3$, $|k_2-m_2|\geq 3$. Applying Lemma \ref{LemD.5}, we get

\begin{equation}\label{P.19}
 \  \| \|Q_{k_1} P_{k_{2}}(D)f \|_{L^2}\|_{\ell_{k_2}^{2,s}
 \ell_{k_1}^{q,a}}
\  \leq C \| \|\widetilde {P}_{m_{2}}(D)\widetilde {Q}_{m_1} f
\|_{L^2}\|_{\ell_{m_2}^{2,s} \ell_{m_1}^{q,a} }.
\end{equation}
For $ 1\leq q \leq 2$, we have the inequality
$$
\| \|\widetilde {P}_{m_{2}}(D)\widetilde {Q}_{m_1} f
\|_{L^2}\|_{\ell_{m_2}^{2,s} \ell_{m_1}^{q,a} } \leq \|
\|\widetilde {P}_{m_{2}}(D)\widetilde {Q}_{m_1} f
\|_{L^2}\|_{\ell_{m_1}^{q,a} \ell_{m_2}^{2,s}  }$$ and from
relation
$$
\| \|\widetilde {P}_{m_{2}}(D)g \|_{L^2}\|_{ \ell_{m_2}^{2,s}  }
\cong \|g\|_{\dot{H}^s}$$ we get
\begin{equation}\label{P.20}
\begin{cases}
\ \| \|Q_{k_1} P_{k_{2}}(D)f \|_{L^2}\|_{\ell_{k_2}^{2,s}
 \ell_{k_1}^{q,a}} \leq \\
\  \leq C \| \|\widetilde {P}_{m_{2}}(D)\widetilde {Q}_{m_1} f
\|_{L^2}\|_{\ell_{m_1}^{q,a} \ell_{m_2}^{2,s}  }\cong\\
\ C \ \| \||D|^{s}Q_{m_1} f \|_{L^2}\|\|_{\ell_{m_1}^{q,a}}\cong
\| f\|_{\ell_x^{q,a}\dot {H}^s}.
\end{cases}
\end{equation}
This inequalities and \eqref{P.9}  imply
\begin{equation}\label{P.21}
\| f\|_{l_D^{2,0}l_x^{q,\alpha}\dot {H}^s}\leq C \| f\|_{l_x^{q,\alpha}\dot {H}^s}.
\end{equation}
This completes the proof. \hfill$\Box$

Further, we obtain in a similar way the following.
\begin{lem}\label{LemP.1dual}
If $q\in[2,\infty]$ and $a , s \in \R$ satisfy
\begin{equation}\label{P.7dual}
\begin{cases}
\  |s|\leq 1,  \\
 \  |a|+|s|<\frac{n}{2}
\end{cases}
\end{equation}
then
\begin{equation}\label{P.8dual}
\| f\|_{\ell_x^{q,a} \dot {H}^s} \leq C\|
f\|_{\ell_D^{2,0}\ell_x^{q,a}\dot {H}^s}
\end{equation}

\end{lem}

In a similar way we can verify the following.

\begin{lem}\label{LemP.1ps}
If $q\in[1,2]$ and $a , s \in \R$ satisfy
\begin{equation}\label{P.7bis}
\begin{cases}
\  |s|\leq1,  \\
 \  |a|+|s|<\frac{n}{2}
\end{cases}
\end{equation}
and $R$ is a pseudo differential operator with convolution type
symbol homogeneous of degree $0,$ then
\begin{equation}\label{P.8bis}
\|R f\|_{l_x^{q,a}\dot {H}^s}\leq C \| f\|_{l_x^{q,a}\dot {H}^s}
\end{equation}

\end{lem}

By using a duality argument one can relax the assumptions on $q$
and obtain the following.

\begin{lem}\label{LemP.1psin}
If $q\in[1,\infty]$ and $a , s \in \R$ satisfy
\begin{equation}\label{P.7bisb}
\begin{cases}
\  |s|\leq1,  \\
 \  |a|+|s|<\frac{n}{2}
\end{cases}
\end{equation}
and $R$ is a pseudo differential operator with convolution type
symbol homogeneous of degree $0,$ then
\begin{equation}\label{P.8bisb}
\|R f\|_{\ell_x^{q,a}\dot {H}^s}\leq C \| f\|_{\ell_x^{q,a}\dot
{H}^s} .
\end{equation}

\end{lem}

\section{Appendix: The Kenig,Ponce, Vega estimate \eqref{eq.sm} for the free Schr\"odinger equation. }

\label{proofKPV}

In this section  we shall recall the basic scale invariant
smoothing estimate due to Kenig,Ponce,Vega.

One possible  proof of the Kenig,Ponce,Vega estimate \eqref{eq.sm}
 is based on the following lemmas:

\begin{lem}\label{Ho1}
For any $u \in  S(\mathbb{R}^n)$ we have
\begin{equation}\label{eq.inclusion}
    \| u(x_1, x^{\prime}) \|_{L^2_{x_1} L^2_{x^{\prime}}}  \leq C \left( \sum_{k \in
\mathbb{Z}} \||x|_k^{1/2} u(x) \|_{L^2_x} \right) ,
\end{equation}
where  $ x = (x_1, x^{\prime}) $ with $x_1\in \mathbb{R}$ and
$x^{\prime}\in \mathbb{R}^{n-1}$.
\end{lem}
\begin{proof} We can consider the case of $n=1$, since
a similar argument works for $n>1$. Let  $u \in S(\mathbb{R})$, we
have
\begin{equation}\label{eq.inclusiona}
    \| u(x) \|_{L^1}  \leq C \left \|  \sum_{k \in
\mathbb{Z}} Q_k(x) u(x)  \right\|_{L^2}.
\end{equation}
From the Cauchy-Schwartz inequality and  the fact that for  the
functions  $Q_k(x)$ we  have $ {\rm supp}_s Q_k(x)\subset
\{2^{k-1}\leq |x|\leq 2^{k+1}\}$, we obtain
\begin{equation}\label{eq.inclusionb}
    \left \|  \sum_{k \in
\mathbb{Z}} Q_k(x) u(x)  \right\|_{L^2}   \leq C \left( \sum_{k
\in \mathbb{Z}} \||x|_k^{1/2} u(x) \|_{L^2} \right) ,
\end{equation} so
\begin{equation}\label{eq.inclusionc}
 \| u(x) \|_{L^1}    \leq C \left( \sum_{k \in
\mathbb{Z}} \||x|_k^{1/2} u(x) \|_{L^2} \right) .
\end{equation}
\end{proof}

Similarly, we have
\begin{lem}\label{Ho2}
For any $u \in  S(\mathbb{R}^n)$ and any $\mathbb{V}\in\R^n$ we
have
\begin{equation}\label{eq.inclusion2}
    \sup_{k \in \mathbb{Z}}  \||x|_k^{-1/2} u(x) \|_{L^2_x}  \leq C \| u(x_1, x^{\prime}) \|_{L^{\infty}_{x_1} L^2_{x^{\prime}}},
\end{equation}
where  $ x = (x_1, x^{\prime}) $ with $x_1\in \mathbb{R}$ and
$x^{\prime}\in \mathbb{R}^{n-1}$.
\end{lem}

The key point in the proof of \eqref{eq.sm} is to establish the
estimate
\begin{equation}\label{eq.strsmoo4}
  \| \partial_1 u(t,x_1, x^{\prime}) \|_{ L^\infty_{x_1} L^2_{t,x^{\prime}}}
  \leq C \| F(t,x_1, x^{\prime}) \|_{ L^1_{x_1} L^2_{t,x^{\prime}}}.
\end{equation}

We shall show now that this estimate completes the proof of
\eqref{eq.sm}.

From \eqref{eq.strsmoo4}, \eqref{eq.inclusion} and
\eqref{eq.inclusion2}  we get
\begin{equation}\label{eq.strsmoo33e}
    \sup_{k \in \mathbb{Z}}  \||x|_k^{-1/2} \partial_1 u(t, x) \|_{L^2_t L^2_x}  \leq
    C \left( \sum_{k \in
\mathbb{Z}} \||x|_k^{1/2} F(t, x) \|_{L^2_t L^2_x} \right) .
\end{equation}
Using the fact that the  Schr\"odinger equation in
\eqref{eq.lwafree}
 and the norm in the right side of
\eqref{eq.strsmoo33e} are invariant under the action of the group
of rotations $ SO(n),$ we obtain
\begin{equation}\label{eq.strsmoo33ej}
    \sup_{k \in \mathbb{Z}}  \||x|_k^{-1/2} \partial_j u(t, x) \|_{L^2_t L^2_x}  \leq
    C \left( \sum_{k \in
\mathbb{Z}} \||x|_k^{1/2} F(t, x) \|_{L^2_t L^2_x} \right) , \ \
\forall j=1,\cdots,n.
\end{equation}

To prove \eqref{eq.strsmoo4} we need of the following two lemmas.

\begin{lem}\label{resolvent1}
There exists a constant $C\geq 0$ so that for all $\lambda \in
\mathbb{C}$ and all $v\in S(\mathbb{R})$ we have
\begin{equation}\label{eq.resolvent1}
  \| v(x) \|_{ L^\infty} \leq C \left \| \left(\frac{d}{dx}-\lambda \right) v(x)\right  \| _{ L^1}.
\end{equation}
\end{lem}

\begin{proof}
Suppose that $Re \lambda\leq 0$. Let $w(x)=(d/dx-\lambda)v(x) \in
S(\mathbb{R}). $ Then
\begin{equation}\label{eq.solution}
  v(x)=\int _{-\infty}^{x}e^{\lambda(x-y)}w(y)dy.
\end{equation}
Thus
\begin{equation}\label{eq.solutiona}
 \| v(x) \|_{ L^\infty} \leq C \left \|w(x)\right  \| _{ L^1}.
\end{equation}

The estimate \eqref{eq.resolvent1} follows for $Re \lambda\leq 0$.
A similar argument works for $Re \lambda\geq 0$. In fact, we have
 $w(x)=(d/dx-\lambda)v(x)$ so
\begin{equation}\label{eq.solution2}
  v(x)=\int ^{\infty}_{x}e^{\lambda(x-y)}w(y)dy.
\end{equation}
and we obtain, as in the previous case,
\begin{equation}\label{eq.solutiona2}
 \| v(x) \|_{ L^\infty} \leq C \left \|w(x)\right  \| _{ L^1}.
\end{equation}

\end{proof}

The estimate \eqref{eq.strsmoo4} follows from the following.

\begin{lem}\label{resolvent2}
For $n \geq 1$ there exists a constant $C\geq 0$ so that for all
$\lambda \in \mathbb{C}$ and all $v\in S(\mathbb{R}^n)$ we have
\begin{equation}\label{eq.resolvent2}
  \| \partial_1 v(x_1, x^{\prime}) \|_{ L^\infty_{x_1}L^2_{ x^{\prime}}}
  \leq C \left \| \left(-\Delta-\lambda \right) v(x_1, x^{\prime})\right  \| _{ L^1_{x_1}L^2_{ x^{\prime}}},
\end{equation}
where  $ x = (x_1, x^{\prime}) $  and
$x^{\prime}=(x_2,\cdots,x_n)\in \mathbb{R}^{n-1}$.
\end{lem}
\begin{proof}
Consider first the case $n=1$. Let be $\lambda=-\mu^2$. Then by
lemma \ref{resolvent1},
\begin{eqnarray}\label{eq.smoo1}
\left  \|  \frac{d}{dx} v(x) \right  \|_{ L^\infty} &\leq&
\frac{1}{2} \left \| \left(\frac{d}{dx}-\mu \right) v(x)\right  \|
_{ L^\infty}+\frac{1}{2} \left \|  \left(\frac{d}{dx}+\mu \right)
v(x) \right \| _{ L^\infty}\leq
\nonumber\\
&\leq& C \left \| \left(-\frac{d^2}{dx^2}-\lambda \right)
v(x)\right  \| _{ L^1}.
\end{eqnarray}
Consider now the  case $n >1.$ Given any $v\in S(\mathbb{R}^n)$ we
denote by $$\widetilde{v}(x_1, k^\prime), \ k^\prime =
(k_2,...,k_n)$$ its partial Fourier transform with respect to
$x^\prime$, i.e.
\begin{equation}\label{eq.fourier}
\widetilde{v}(x_1, k^\prime)=(2\pi)^{-(n-1)}\int e^{-i k^\prime
x^\prime} v(x_1,x^\prime) dx^\prime.
\end{equation}
Using the one-dimensional result \eqref{eq.smoo1}, for each fixed
$k^\prime$, we obtain
\begin{equation}\label{eq.fourier2}
|\partial_1 \widetilde{v}(x_1, k^\prime)|^2 \leq C \int
\left|(-\partial^2_1+|k^\prime|^2 - \lambda) \widetilde{v}
\right|^2 dx_1.
\end{equation}
Integrating with respect $ k^\prime$ and using the Plancherel
identity, we derive
\begin{equation}\label{eq.resolvent3}
  \| \partial_1 v(x_1, x^{\prime}) \|_{ L^\infty_{x_1}L^2_{ x^{\prime}}}
  \leq C \left \| \left(-\Delta-\lambda \right) v(x_1, x^{\prime})\right
   \| _{ L^1_{x_1}L^2_{ x^{\prime}}}.
\end{equation}
This completes the proof.
\end{proof}

In fact the basic idea of the proof of \eqref{eq.strsmoo4} is to compute the Fourier transform
with respect the temporal variable of the equation \eqref{eq.lwafree} and obtain
\begin{equation}\label{eq.resfin}
 - i\Delta \hat{u}(\lambda,x) - i\lambda \hat{u}(\lambda,x) =\widehat{F}(\lambda,x) .
\end{equation}
Now if we split $ x = (x_1, x^{\prime})\in  \mathbb{R}\times \mathbb{R}^{n-1}$
and indicate with $v_{\lambda}(x_1, x^{\prime})=\hat{u}(\lambda,x)$,
we can apply the Lemma \ref{resolvent2} and have

\begin{eqnarray}\label{eq.resolventfin2}
  \| \partial_1 v_{\lambda}(x_1, x^{\prime}) \|_{ L^\infty_{x_1}L^2_{ x^{\prime}}}
  \leq C \left \| \left(-\Delta-\lambda \right) v_{\lambda}(x_1, x^{\prime})\right
   \| _{ L^1_{x_1}L^2_{ x^{\prime}}}.
\end{eqnarray}
This estimate with the equation \eqref{eq.resfin} give the following other one

\begin{eqnarray}\label{eq.resolventfin3}
  \| \partial_1 \hat{u}(\lambda, x_1, x^{\prime}) \|_{ L^\infty_{x_1}L^2_{ x^{\prime}}}
  &\leq& C \left \| \left(-\Delta-\lambda \right) \hat{u}(\lambda, x_1, x^{\prime})\right
   \| _{ L^1_{x_1}L^2_{ x^{\prime}}}\leq
   \nonumber\\
  \\
  &\leq&  \left \| \widehat{F}(\lambda, x_1, x^{\prime})\right
   \| _{ L^1_{x_1}L^2_{ x^{\prime}}}.
   \nonumber
\end{eqnarray}
The application of Plancherel's theorem in the temporal variable
 gives the estimate \eqref{eq.strsmoo4} and this completes the proof of \eqref{eq.sm}.

\bibliographystyle{plain}

\end{document}